\documentstyle[12pt]{article}

\title{Self-Dual Manifolds with\\ Positive Ricci Curvature}

\author{
Claude LeBrun\thanks{Supported in part by NSF grant DMS-9204093}\\
{ Department of Mathematics} \\
{ State University of New York}\\
{ Stony Brook, NY 11794, USA}
\bigskip\\
Shin Nayatani\\
{ Max-Planck-Institut f\"{u}r Mathematik}\\
{ Gottfried-Claren-Stra{\ss}e 26} \\ { D-53225 Bonn, Germany}\\
{ and}\\
{ Mathematical Institute}\\ {T\^{o}hoku University}\\
{ Sendai 980, Japan}
\bigskip\\
Takashi Nitta\\
{ Department of Mathematics} \\ { Mie University}\\
{ Tsu 514, Japan}
}

\date{}

\newtheorem{thm}{Theorem}[section]
\newtheorem{prop}[thm]{Proposition}
\newtheorem{Lemma}[thm]{Lemma}
\newtheorem{cor}[thm]{Corollary}
\newtheorem{defn}{Definition}[section]
\def\hook{\mbox{}\begin{picture}(10,10)\put(1,0){\line(1,0){7}}
  \put(8,0){\line(0,1){7}}\end{picture}\mbox{}}
\def\laph{\triangle}
\def\lapg{\Delta}
\def\schou{{\rm Q}}
\def\e{{\bf e}}
\def\cpt{M}
\def\bdl{{\cal M}}
\def\disk{D^3}

\newcommand{\Ric}{\mathop{\rm Ric}\nolimits}

\newenvironment{proof}{\medskip
\noindent {\bf Proof.}}{\hfill \rule{.5em}{1em}\mbox{}\bigskip}

\newcounter{remark}[section]
\renewcommand{\theremark}{\thesection.\arabic{remark}}

\newcommand{\sctn}{\setcounter{equation}{0}
 \setcounter{remark}{0}\section}

\begin{document}

\maketitle

\begin{abstract}
 We prove that the connected sums
${\bf CP}_2\# {\bf CP}_2$ and ${\bf CP}_2\# {\bf CP}_2\# {\bf CP}_2$
 admit self-dual metrics
with positive Ricci curvature.
Moreover, every self-dual metric of positive scalar curvature
on  ${\bf CP}_2 \#{\bf CP}_2$  is conformal to a metric
with positive Ricci curvature.
\end{abstract}

\newpage

\section*{Introduction.}

An oriented Riemannian 4-manifold  $(M,g)$ is said to be
{\em self-dual} if its Weyl
curvature  $W$, thought of as a bundle-valued 2-form,
satisfies $W=\star W$, where  $\star$ denotes the
 Hodge star operator. Because both $W$ and $\star$ are
unchanged if the metric is multiplied by a
positive function, this property is  conformally
invariant,
 and  the term self-dual is thus  often
used to describe the  conformal class
$[g]:= \left\{ ug ~|~ u: M \stackrel{C^{\infty}}{\longrightarrow}
{\bf R}^+\right\}$ rather than  the
 metric $g$ which represents it.

Two  familiar
 examples of compact self-dual manifolds are
provided by the  symmetric spaces
$S^4=SO(5)/SO(4)$ and ${\bf CP}_2=SU(3)/U(2)$.
For many years, these were the only known examples
of compact simply-connected
 self-dual manifolds with positive scalar curvature,
and it was therefore a major breakthrough when
Poon \cite{poo}  constructed a one-parameter family of
positive-scalar-curvature
self-dual metrics
on  ${\bf CP}_2 \#{\bf CP}_2$;
here the connected sum operation $\#$
is carried out by deleting balls from the given
manifolds and then identifying the resulting
boundaries in a manner compatible with the given
orientations.
Motivated by this discovery, Donaldson-Friedman \cite{dofr}
and Floer \cite{flo} abstractly
 proved the existence of self-dual metrics on
the $n$-fold connected sum
$$n{\bf CP}_2:= \underbrace{{\bf CP}_2 \#\cdots \#{\bf CP}_2}_n$$
for every $n$.
The first author \cite{leb2} then
realized that  such metrics on $n{\bf CP}_2$
can be  constructed
explicitly
 by means of the so-called ``hyperbolic ansatz'' reviewed below in
 \S \ref{ansatz}. This last method  has the added
advantage that each of the   conformal classes
so constructed  can
be seen to contain a representative of  positive scalar curvature.

On the other hand,  the
examples provided by
 $S^4$ and ${\bf CP}_2$ actually have positive {\em Ricci}
curvature, and, in light of the work of
Cheeger \cite{ch}, Anderson \cite{and} and Sha-Yang \cite {sy},
 it is natural to ask whether there are
other compact 4-manifolds which admit self-dual
metrics with this property. Our objective here is to show  that the
answer to this question is {\em yes}. We will accomplish this
(\S 4) by explicitly  constructing  such metrics on
 $n{\bf CP}_2$ when $n=2$ and $3$; moreover,
  it will turn out that
each of Poon's  conformal classes on
$2{\bf CP}_2$   contains such a metric.
On the other hand, it is rather
easy (\S 1) to see that a compact self-dual
manifold with positive Ricci curvature must be
 homeomorphic to $n{\bf CP}_2$ for some $n \geq 0$,
where by convention   $0{\bf CP}_2 := S^4$.
This raises the
fascinating question, left unanswered here,
 of whether
$n{\bf CP}_2$ admits such  metrics when $n \geq 4$.

On a related  front, Gauduchon \cite{gau} has studied self-dual
manifolds with  {\sl non-negative Ricci operator} (cf. \S 1), and asked
whether $2{\bf CP}_2$ and $3{\bf CP}_2$ admit such metrics.
Our metrics on $2{\bf CP}_2$ will be seen to  satisfy both
this condition and another,  which we call
{\sl strongly positive Ricci curvature}.

In order to prove these positivity results,
we will first  (\S 2) need to  compute the  Ricci curvature
of the general  self-dual metric of  hyperbolic-ansatz type.
Our results will  then follow once
we have introduced  a suitable choice of conformal
gauge, motivated (\S 3) by a re-examination  of
 the Fubini-Study metric of
${\bf CP}_2$.

\bigskip

\noindent
{\bf Acknowledgements.} The
present article grew out of a series of conversations
at the Mathematical Sciences Research Institute
and the Mathematics Department of UC Berkeley; we are thus greatly
 indebted to these institutions for their  hospitality and
financial support. The  first and second author would
also  like to respectively thank the Erwin Schr\"odinger Institute
(Vienna) and the
 Max-Planck-Institut f\"{u}r
Mathematik (Bonn) for their financial
 support and  hospitality.
Finally, the third author would like
to  thank  Prof. S. Kobayashi
for sponsoring his visit to Berkeley.

\bigskip

\sctn{Topological Preliminaries}

The  present article is largely motivated by
 the following easy observation:

\begin{prop}\label{pr1.1}
Let $(M, g)$ be a compact self-dual 4-manifold with positive Ricci
curvature.
Then $M$ is homeomorphic to $n{\bf CP}_2$ for some $n \geq 0$.
Moreover,  $M$ is  {\em diffeomorphic} to $n{\bf CP}_2$
if $n\leq 4$.
\end{prop}

\begin{proof}  Let  the universal cover $\widetilde{M}$
of $M$ be equipped with the pull-back metric.
Since  the Ricci curvature of $\widetilde{M}$ is
then bounded below by
a positive constant,
 Myers' theorem tells us that $\widetilde{M}$ is compact,
 and $\widetilde{M}\to M$ is therefore a finite-sheeted covering.
However, a simple Bochner-Weitzenb\"ock argument  \cite{bou,leb}
implies  that a compact self-dual 4-manifold with positive
scalar curvature must have  $b_- = 0$.
Thus for both $M$ and $\widetilde{M}$, we have $b_1 = b_- = 0$,
 and hence both have
$\chi - \tau = 2(1 - b_1 + b_-) = 2$, where $\chi$ is the Euler
characteristic and $\tau$ is the signature.
But $\chi - \tau$ is multiplicative under finite coverings because
it can be computed from a Gauss-Bonnet formula.
Hence  $\widetilde{M}\to M$ is the trivial  covering,
 and $M$ is simply connected.
It now follows from the work of Donaldson
\cite{D} and Freedman
 that $M$ is
homeomorphic to $n{\bf CP}_2$ for some $n \geq 0$.
On the other hand, a self-dual manifold with
positive scalar curvature, $b_1=0$, and $\tau \leq 4$
must \cite{pepo} be  {\em  diffeomorphic}   to
$n{\bf CP}_2$ because its twistor space
contains a rational hypersurface of degree 2.
\end{proof}
\bigskip

It is now natural to ask whether, conversely, the manifolds
 $n{\bf CP}_2$ admit self-dual metrics with positive
Ricci curvature. For small values of $n$ we shall see that the
answer is in fact affirmative.
\bigskip

Rather than merely asking for  the Ricci curvature $\Ric$
to   be positive, one might ask for its trace-free part
$\Ric_0$ to be  small enough with respect to
its scalar curvature $s > 0$ so as to   guarantee {\em a priori}
 that $\Ric > 0$. This motivates the following definition:

\begin{defn}
Let $(M,g)$ be a Riemannian 4-manifold. Then we will
say that $M$ has {\em strongly positive Ricci curvature}
if, at each point of $M$, we have
$$| \Ric_0 | < \frac{s}{2\sqrt{3}}.$$
Similarly, we will say that $M$ has
{\em strongly non-negative Ricci curvature} if $s >0$
and
$$| \Ric_0 | \leq \frac{s}{2\sqrt{3}}$$
at every point of $M$.
\end{defn}

Observe that strongly positive Ricci curvature implies
  positive Ricci curvature. Indeed,
if $(\lambda_1, \ldots , \lambda_4)$ are
the eigenvalues of $\Ric / s$, then, in the 3-plane
$\lambda_1 + \cdots + \lambda_4 =1$, positive
Ricci curvature
corresponds to the tetrahedron with corners
$(1,0,0,0), \ldots , (0,0,0,1)$, whereas
strongly positive Ricci curvature corresponds to the
ball of radius $1/2\sqrt{3}$ around  $(\frac{1}{4},
\ldots, \frac{1}{4})$, and  this ball  just fills the
 in-sphere of the tetrahedron.  By the
same argument, we also see that strongly
non-negative Ricci curvature implies
non-negative Ricci curvature.

\begin{prop}
Let $(M,g)$ be a self-dual 4-manifold with strongly
positive Ricci curvature. Then $M$ is diffeomorphic
to $n{\bf CP}_2$, where $0\leq n\leq 3$.
\end{prop}

\begin{proof}
The Gauss-Bonnet formulae for the signature
and Euler characteristic of a compact  oriented Riemannian
4-manifold $(M, g)$ are
$$\tau (M)
= \frac{1}{12\pi^2}\int_M \left( |W_+|^2- |W_-|^2\right) v_g$$
and
$$\chi (M) = \frac{1}{8\pi^2}\int_M \left(
|W_+|^2+|W_-|^2-\frac{|\Ric_0|^2}{2}+ \frac{s^2}{24}\right) v_g ,$$
where $v_g$ is the metric volume form.
Thus any compact self-dual 4-manifold satisfies
$$
(2\chi  - 3\tau ) (M)= \frac{1}{8\pi^2} \int_M \left(
\frac{s^2}{12}-{|\Ric_0|^2}\right) v_g ,$$
and the right-hand side  is manifestly
positive if the Ricci curvature is strongly positive.
Now $M$ is homeomorphic to $n{\bf CP}_2$ by
Proposition \ref{pr1.1}, and even diffeomorphic if
$n\leq 4$.
But the inequality
$2\chi - 3\tau > 0$ implies that $n < 4$, as desired.
\end{proof}

Rather than focusing on the Ricci tensor of
 a Riemannian $4$-manifold $(M,g)$, one
may  instead choose \cite{gau} to consider
 an algebraically equivalent
object ${\cal R}ic$, called the {\em Ricci operator},
which is  defined as the full curvature operator
minus its Weyl component.
If we let $\schou$ denote Schouten's  modified Ricci tensor
$$
\schou =   \Ric - \frac{s}{6}  g ,
$$
the Ricci operator is explicitly given by
$$
{\cal R}ic (X \wedge Y) = \frac{1}{2}
\left( \schou^\sharp (X) \wedge Y + X \wedge \schou^\sharp (Y) \right),
$$
where  $\schou^\sharp$
is the endomorphism of
$TM$ corresponding to $\schou$
and $X$, $Y$ are any tangent vectors.
It follows that the Ricci operator
is positive (respectively,  non-negative)
if and only if the sum of the lowest two eigenvalues of $\schou$
is positive (respectively,  non-negative).
In terms of  $\lambda_1, \ldots
, \lambda_4$,  this  corresponds
  to   requiring that
$$\frac{2}{3} >
(\mbox{resp. } \geq )~ \lambda_i + \lambda_j >
(\mbox{resp. } \geq )~  \frac{1}{3} \quad
\forall i\neq j, $$
which is to say that
$(\lambda_1, \cdots, \lambda_4)$ is a point of
the open (respectively, closed) cube with
corners $(\frac{1}{3},  \frac{1}{3},\frac{1}{3},0),
\ldots , (0, \frac{1}{3},\frac{1}{3},\frac{1}{3}),
(\frac{1}{6},\frac{1}{6},\frac{1}{6},\frac{1}{2}),
\ldots ,  (\frac{1}{2}, \frac{1}{6},\frac{1}{6},\frac{1}{6})$.
Since this cube is contained in the
in-sphere, we therefore have
\begin{quote}
positive Ricci operator $\Rightarrow$\\ strongly positive
Ricci curvature
$\Rightarrow$\\  positive
Ricci curvature,
\end{quote}
and
\begin{quote}
non-negative Ricci operator and $s> 0$
$\Rightarrow$ \\strongly non-negative
Ricci curvature
$\Rightarrow$ \\ non-negative
Ricci curvature.
\end{quote}
Moreover, non-negative Ricci operator and $s> 0$
fail to imply that  the
Ricci curvature is strongly positive only when
$(\lambda_1, \ldots , \lambda_4)$ is a corner of the cube.
Using this observation, we  now prove a slightly sharpened version
of a result discovered by Gauduchon \cite{gau}, using
different  methods.

\begin{thm}\label{th1.2}
Let $(M, g)$ be a compact self-dual 4-manifold with positive
scalar curvature and non-negative Ricci operator.
Then either
 $M$ is
diffeomorphic to
$n{\bf CP}_2$, $0 \leq n \leq 3$, or else the
universal cover of $(M, g)$ is
 the Riemannian product ${\bf R} \times S^3$.
\end{thm}

\begin{proof}
Since the Ricci curvature is strongly non-negative,
$$(2\chi  - 3\tau ) (M)= \frac{1}{8\pi^2} \int_M \left(
\frac{s^2}{12}-{|\Ric_0|^2}\right) v_g \geq 0,$$
with equality iff $|\Ric_0 | \equiv s/ 2\sqrt{3}$.
If the inequality is strict,
  $2\chi - 3 \tau > 0$.
Thus $b_1(M)= 0$ and $\tau (M) < 4$.
The proof of Proposition \ref{pr1.1} thus implies that
$M\approx n{\bf CP}_2$ for $n < 4$.

If equality holds, $(\lambda_1, \ldots , \lambda_4)$
must everywhere be one of the corners of the
previously mentioned cube, and
$\Ric$ therefore has exactly two eigenvalues at each
point of $M$, one with
multiplicity 3 and one with multiplicity 1. It follows that
there is a line sub-bundle of $TM$, and
$\chi (M)= 0$. Moreover,
  $b_+(M)= \tau (M) = \frac{3}{2} \chi (M) =0$, so that
 $b_2(M)= b_+(M) = 0$.
Hence $0=\chi (M) = 2 - 2b_1(M)$, and $b_1(M)= 1$.
Since $M$ has non-negative Ricci curvature, the classical
Bochner argument \cite{bes} now says that
$M$ admits a parallel 1-form, and thus locally splits as
the Riemannian product of ${\bf R}\times N$, where $N$ is
 a 3-manifold. But since $\Ric$ everywhere has a
 positive eigenvalue of multiplicity 3,
$N$  is an Einstein 3-manifold of positive scalar curvature.
Thus $N$ has positive constant sectional
curvature, and the universal cover of  $M$ is ${\bf R}\times S^3$.
\end{proof}

\sctn{Ricci Curvature and the Hyperbolic Ansatz}\label{ansatz}

In this section, we shall compute the Ricci curvature of
those self-dual metrics which arise from the following
``hyperbolic ansatz'' construction:

\begin{prop}\label{pr2.1} {\rm \cite{leb2}}
Let $({\cal H}^3, h)$ denote hyperbolic 3-space, which we equip with
 a fixed
orientation, and let $V$ be a
positive harmonic function on some
 open set ${\cal V} \subset {\cal H}^3$.
Suppose that the cohomology class of $\frac{1}{2\pi} \star dV$
is integral, where $\star$ is the Hodge star operator of ${\cal H}^3$.
Let $\bdl \rightarrow {\cal V}$ be a circle bundle with a connection
1-form $\theta$ whose curvature is $\star dV$.
Then the conformal class
$$
[g] = \left[Vh + V^{-1} \theta^2\right]
$$
of Riemannian metrics on $\bdl$
is  self-dual with respect to the orientation
 determined by $\theta \wedge  v_h$, where
$v_h$ is the volume form of ${\cal H}^3$.
\end{prop}

We now wish to calculate the Ricci curvature of
metrics in these self-dual conformal classes.
With the most obvious choice of conformal factor,
the answer turns out to be surprisingly simple:

\begin{prop} \label{naive}
For any positive harmonic function $V$ on a region of ${\cal H}^3$,
the
Ricci curvature of the self-dual metric $g=Vh+V^{-1}\theta^2$
is  $\Ric_g=-2h$.
\end{prop}
The $V$-independence of this Ricci curvature
 is analogous to the
Ricci-flatness of the  metrics produced
via the Gibbons-Hawking
ansatz \cite{giha}.

While this answer is beguilingly simple, it is also
depressingly negative!
Fortunately, the picture will become  less bleak
once we   conformally
rescale  our  metric:

\begin{prop}\label{general}
Let $f$ and $V$ be respectively a smooth function and
a positive harmonic function  on a
domain ${\cal V}\subset {\cal H}^3$. Then
the Ricci curvature of the corresponding self-dual metric
$g = e^{2f} \left(V h + V^{-1} \theta^2\right)$ is given by
\begin{eqnarray}
\Ric_g &=& \left(-2 - \laph f - 2|df|^2 - V^{-1} \langle
dV, df\rangle \right) h
\nonumber\\
&& - 2D df + 2(df)^2 + 2V^{-1} dV \odot df \label{eq12}\\
&& + \left(-\laph f - 2|df|^2 + V^{-1} \langle dV, df\rangle \right)
(V^{-1}\theta)^2
\nonumber\\
&& - 2V^{-1} \star (dV \wedge df) \odot V^{-1}\theta. \nonumber
\end{eqnarray}
Here $D$, $\laph$, and $\star$
are respectively
 the Levi-Civit\`a connection, negative
Laplace-Beltrami operator, and Hodge star operator
of   hyperbolic 3-space $({\cal H}^3, h)$, while
 $| \cdot |$ and $\langle \cdot , \cdot \rangle$
are  the corresponding norm and inner product
on 1-forms.
\end{prop}

To prove these statements, let us
first observe that (\ref{eq12}) is valid iff it holds for
some {\em particular} $f$; in particular, Propositions
\ref{naive} and \ref{general} are logically equivalent.
Indeed, if $g_0= Vh+V^{-1}\theta^2$ and
$g=e^{2f}g_0$, the standard formula \cite{bes}
governing the alteration of curvature by
conformal rescaling yields
$$\Ric_g=\Ric_{g_0}- 2\nabla df + 2(df)^2
-\left(\lapg f + 2|df|_{g_0}^2\right)g_0 ,$$
 where $\nabla$ and $\lapg$ are respectively
the Levi-Civit\`a connection and negative
Laplace-Beltrami operator of $g_0$.
Now since
\begin{eqnarray*}
\nabla df&=& {\textstyle \frac{1}{2}}
\pounds_{\mbox{\tiny grad}_{g_0}f}g_0 \\
&=& {\textstyle \frac{1}{2}}
\pounds_{\mbox{\tiny grad}_{g_0}f}Vh
+{\textstyle \frac{1}{2}}
\pounds_{\mbox{\tiny grad}_{g_0}f}V^{-1}\theta^2\\
&=& {\textstyle \frac{1}{2}}
\pounds_{V^{-1}\mbox{\tiny grad}_hf}Vh
+ {\textstyle \frac{1}{2}}
\left(\pounds_{V^{-1}\mbox{\tiny grad}_hf}V^{-1}\right)\theta^2
+ V^{-1} \theta \odot \left(V^{-1}\mbox{grad}_hf \hook d\theta\right)\\
&=&  {\textstyle \frac{1}{2}}
\left(\pounds_{V^{-1}\mbox{\tiny grad}_hf}V\right)h
+ V \mbox{symm}(D V^{-1} df)
 -  \frac{\langle dV , df \rangle}{2V^3} \theta^2
+ \frac{\theta \odot (\mbox{grad}_hf \hook \star  dV)}{V^2}\\
&=& Ddf - \frac{dV\odot df}{V}
+ \frac{\langle dV , df \rangle}{2V}  h
 - \frac{\langle dV , df \rangle}{2V^3} \theta^2
+\frac{\theta \odot \star (dV \wedge df )}{V^2}, \end{eqnarray*}
it follows that
$$ \lapg f = V^{-1} \laph f ,$$
and we therefore have
\begin{eqnarray*}
\Ric_g&=&\Ric_{g_0}- 2  Ddf +2V^{-1} dV\odot df
- V^{-1} \langle dV, df\rangle  h\\&&
 +  V^{-3} \langle dV , df \rangle \theta^2
-2V^{-2} \theta \odot \star (dV \wedge df ) + 2(df)^2 \\&&
-\left(V^{-1}\laph f + 2V^{-1}|df|^2\right)(Vh+V^{-1}\theta^2)
\\&=& \Ric_{g_0} - \left(\laph f + 2|df|^2+V^{-1}
 \langle dV, df\rangle\right)h\\
&&- 2  Ddf + 2(df)^2+2V^{-1} dV\odot df\\&&
+ \left(- \laph f- 2|df|^2+ V^{-1}
\langle dV , df \rangle\right)(V^{-1}\theta)^2
\\&& -2V^{-2} \theta \odot \star (dV \wedge df ).
\end{eqnarray*}
But this will coincide with (\ref{eq12}) for any  particular $f$
iff  $\Ric_{g_0}=-2h$.

We now complete our proof by
verifying (\ref{eq12}) for  a slightly peculiar choice of $f$,
 best described in terms of the
 upper-half-space model
$$h=\frac{dx^2+dy^2+dz^2}{z^2}, ~~~~ z > 0,$$
of ${\cal H}^3$. We will now set $f=\log z$
because \cite[\S 3]{lhopf} the corresponding metric
$$g=z^2(Vh+V^{-1}\theta^2)$$
is K\"ahler
with respect to the integrable almost-complex structure
$$dx\mapsto dy, ~~~~ dz\mapsto \frac{z}{V}\theta,$$
with
Ricci form
$$P=-d(V^{-1}\theta)=  -\frac{\star dV}{V}+
\frac{dV\wedge \theta}{V^2}.$$
The Ricci curvature of this metric is therefore
\begin{eqnarray*}
\Ric_g &=& \frac{V_z}{zV}\left[-dx^2-dy^2+dz^2
+\left(\frac{z}{V}\theta \right)^2\right]
+ \frac{2V_x}{zV} \left[ dx\odot dz +  dy \odot \frac{z}{V}\theta\right]
\\&&
+ \frac{2V_y}{zV} \left[dy\odot dz -  dx \odot\frac{z}{V} \theta \right].
\end{eqnarray*}
But, since $|df|^2= 1$,
$$Ddf =  \frac{1}{2}\pounds_{\mbox{\tiny grad}_h f} h=
\frac{1}{2}\pounds_{z\frac{\partial}{\partial z}}
 \left(\frac{dx^2+dy^2+dz^2}{z^2}\right)=-\frac{dx^2+dy^2}{z^2},$$
and $\laph f = -2$,
this is exactly the result predicted by
(\ref{eq12}) with $f=\log z$. Thus (\ref{eq12}) holds
for our particular $f$, and
 Propositions \ref{naive} and \ref{general} therefore follow.

To conclude this section, let us
point  out that the scalar curvature $s_g$ and the modified
Ricci  tensor
$\schou_g =   \Ric_g - \frac{1}{6} s_g g$
are now respectively  given by
\begin{equation}\label{eq13}
s_g = 6 e^{-2f} V^{-1} \left(-1 - \laph f - |df|^2\right)
\end{equation}
and
\begin{eqnarray}
\schou_g &=&   \left(-1 - |df|^2 -
\langle \psi , df \rangle \right) h\nonumber\\
&& - 2D df + 2 (df)^2 +   2\psi \odot df\label{hmm}\\
&& +  \left(1 - |df|^2 +   \langle \psi , df \rangle
\right)
(V^{-1}\theta)^2\nonumber\\
&& -   2\star (\psi \wedge df) \odot V^{-1}\theta,\nonumber
\end{eqnarray}
where $\psi = V^{-1}dV= d\log V$.
Notice that  the sign of $s_g$ is independent of
$V$; for applications, cf. \cite{leb2,kim}.

\sctn{Choosing a Conformal Factor}

The hyperbolic ansatz described in the
last section can be used \cite{leb2} to
construct self-dual metrics on $n{\bf CP}_2$.
When $n=1$, this construction gives metrics conformal to the
Fubini-Study metric on ${\bf CP}_2$, and  our main tasks here
will be to re-examine the type of conformal factor
this entails.

Let $\{ p_1, \dots, p_n \}$ be an arbitrary collection of $n$ points
in ${\cal H}^3$, and let
$$
G_j = \frac{1}{2} (\coth r_j -1)
$$
be the hyperbolic Green's function centered at $p_j$;
here $r_j$ is the hyperbolic distance from $p_j$, and our
normalization is chosen so that $d\star d G_j = -2\pi \delta_{p_j}$.
Thus
\begin{equation}
V := 1 +  \sum_{j=1}^n G_j =
1 + \frac{1}{2} \sum_{j=1}^n (\coth r_j - 1)
\label{defV}
\end{equation}
 is a positive harmonic function on
${\cal V} = {\cal H}^3 \setminus \{ p_1, \dots, p_n \}$
satisfying the integrality condition of Proposition \ref{pr2.1}.
Letting $(\bdl, \theta)$ be the circle bundle with connection 1-form
as in Proposition \ref{pr2.1}, which is uniquely determined up to
gauge equivalence since ${\cal V}$ is simply connected, we thus obtain
a self-dual metric
$$
g_0 = Vh + V^{-1} \theta^2
$$
on $\bdl$.
If we now use the Klein projective model
to  identify ${\cal H}^3$ with
the interior of the closed 3-disk
$\disk$, there is a smooth compactification
$\cpt$ of $\bdl$ such that the bundle projection
$\bdl \rightarrow  {\cal H}^3 \setminus \{ p_j \}$
extends to a surjective smooth map
$\cpt \rightarrow \disk$, and $\disk$ is thereby
identified with the orbit space of an $S^1$-action on
$\cpt$; in fact, $\cpt\setminus\bdl$ is the set of fixed points
of this action, and consists of a 2-sphere $\widehat{S}^2$,
which projects diffeomorphically to
  $\partial \disk$, and  $n$ isolated fixed points
  $\widehat{p}_j$, one for each  $p_j\in {\cal H}^3$.
 Moreover,     $g = e^{2f} g_0$ extends
to a self-dual metric on the
compact manifold $\cpt\approx n{\bf CP}_2$ whenever
  $f:{\cal H}^3\to
{\bf R}$ is a smooth function
which behaves like  $-r$ near infinity, where
  $r$ is the hyperbolic distance
from an arbitrary reference point.
When $n = 0, 1$, this construction
produces the  conformal classes of
the standard metrics on $S^4$ and ${\bf CP}_2$;
when $n = 2$, it instead yields
 the self-dual metrics on $2{\bf CP}_2$
first discovered by Poon \cite{poo}.

 In the above discussion, we assumed for simplicity
that $f$ was a smooth function on  ${\cal H}^3$; and on
  ${\cal H}^3\setminus \{p_j\}$ smoothness is obviously
needed to guarantee that $e^{2f}g_0$ is smooth on $\bdl$.
On the other hand, the   derivative of the natural projection
 $\cpt\to \disk$ vanishes at each
$\widehat{p}_j$, and the pull-back of the function
$r_j$ is consequently smooth on $\cpt \setminus S^2$.
Choices of  $f$ with
this sort of behavior near the $p_j$ are also allowable, and
will in fact turn out to be crucial for our purposes.

To see why, let us look more closely at the $n=1$ case.
In geodesic polar coordinates about $p=p_1$, the
hyperbolic metric on
${\cal H}^3\setminus p$ can be
written as
$$h=dr^2+ \sinh^2r~ g_{S^2},$$
where $g_{S^2}$ is the standard metric on the
unit 2-sphere. Now the ansatz  stipulates
that $V= 1+ \frac{1}{2}(\coth r -1) = (1-e^{-2r})^{-1}$,
and hence  $\star dV= -\frac{1}{2}\omega$, where
$\omega$ is the standard area form on the 2-sphere.
 In order to produce a circle bundle with this
curvature, let $\mu : S^3\to S^2$ be the Hopf map,
and let   the unit 3-sphere
 $S^3=Sp(1)$ be equipped with a left-invariant
orthornomal coframe $\{\sigma_1, \sigma_2,\sigma_3\}$
such that
$\mu^*g_{S^2}= 4({\sigma_1}^2+{\sigma_2}^2)$.
Then $\mu^*(-\frac{1}{2}\omega) = -2\sigma_1\wedge\sigma_2=
d(-\sigma_3)$, and the
desired circle bundle $\pi: \bdl\to {\cal H}^3\setminus
p$ may be taken to be  the pull-back of $\mu$,
with connection form $\theta=-\sigma_3$, to
$S^2\times {\bf R}^+$.
Thus
\begin{eqnarray*}
g_0&=&Vh+V^{-1}\theta^2
\\
&=& \frac{1}{1-e^{-2r}}\left[dr^2 + 4\sinh^2 r  ~
( {\sigma_1}^2+ {\sigma_2}^2 ) \right]+ (1-e^{-2r}) {\sigma_3}^2.
\end{eqnarray*}
Setting $\rho= \cos^{-1}(e^{-r})$, we now have
\begin{eqnarray*}
e^{-2r} g_0&=& \cot^2 \rho\left[
\tan^2\rho ~d\rho^2
+ \tan^2\rho~ \sin^2\rho~ ( {\sigma_1}^2+ {\sigma_2}^2 ) \right]+
\cos^2\rho~ \sin^2 \rho ~ {\sigma_3}^2\\&=&
d\rho^2 + \sin^2 \rho ~({\sigma_1}^2+ {\sigma_2}^2
 + \cos^2\rho ~{\sigma_3}^2  ),
\end{eqnarray*}
which is exactly the Fubini-Study metric of
${\bf CP}_2$, expressed in geodesic polar coordinates.
So far as positive Ricci curvature is concerned,
 the  best possible choice of $f$ when $n=1$
is thus $f=-r$, and the  challenge now facing us
is to suitably  generalize this for $n > 1$.  Since
we will still need $f\sim -r$  as $r\to \infty$,
one obvious generalization is
$$f= -\frac{r_1 + \cdots + r_n}{n}.$$
In the next section, we will see that
this choice actually works surprisingly well when  $n\leq 3$.

\sctn{Positive Ricci Curvature}

In the previous section, we associated a conformal class
of self-dual metrics on $n{\bf CP}_2$ to any configuration
of points $\{ p_1, \dots, p_n \}$ in ${\cal H}^3$.
We will henceforth denote this conformal class by
 $C_{p_1, \dots, p_n}$.

\begin{thm}\label{th4.1}
 Each conformal class
$C_{p_1, p_2}$
of self-dual metrics on
${\bf CP}_2\#
{\bf CP}_2$
contains a metric with strongly positive Ricci curvature
and non-negative Ricci operator.
\end{thm}
In fact, the metric $g=e^{2f}(Vh+V^{-1}\theta^2)$ has these
properties provided we set
$$f= -\frac{r_1+r_2}{2} , $$
where $r_1$ and $r_2$ are  respectively the hyperbolic distances from
$p_1,p_2\in {\cal H}^3$.
We will prove this by
 first showing that the Ricci operator is non-negative, and then
observing that the   Ricci curvature is still strongly positive
at the points where the
Ricci operator has non-trivial kernel.

On an open dense subset of $\bdl\subset \cpt$,
and with respect to the metric  $V^{-1}g_0= h+V^{-2}\theta^2$,
we may  define
an oriented orthonormal coframe $\{ \e^1 , \ldots, \e^4\}$
by $$\e^1 = \frac{dr_1+dr_2}{|dr_1+dr_2 |}, ~~
\e^2= \frac{dr_1-dr_2}{|dr_1-dr_2 |}, \mbox{ and }
\e^4= V^{-1}\theta.$$
 Let $\varphi :=\sin^{-1}\langle dr_1, \e^1 \rangle$ be the oriented
angle between $dr_1$ and $\e^1$. Then
\begin{eqnarray*}
df &=& - \left( \cos \varphi \right) \e^1  ,\\&&\\
dV &=& -\frac{1}{2}  \left[
 \frac{dr_1}{\sinh^2 r_1}+\frac{dr_2}{\sinh^2 r_2}\right] \\
&=& -\frac{1}{2} \left[\cos \varphi
\left( \frac{1}{\sinh^2 r_1}+ \frac{1}{\sinh^2 r_2}\right) \e^1 +
\sin \varphi \left( \frac{1}{\sinh^2 r_1}-
\frac{1}{\sinh^2 r_2}\right) \e^2
\right],
 \\&&\\
D df &=& - \frac{1}{2} \left[ \coth r_1 \left(h - d{r_1}^2\right)
+\coth r_2 \left(h - d{r_2}^2\right)\right] \\
&=& - \frac{1}{2}\left[
\sin^2 \varphi \left( \coth r_1 + \coth r_2  \right) (\e^1)^2
-2 \cos \varphi \sin \varphi \left(  \coth r_1-\coth r_2 \right)
\e^1 \odot \e^2\right.\\
&&\left. + \cos^2 \varphi\left( \coth r_1 + \coth r_2 \right)
(\e^2)^2
+ \left(  \coth r_1 + \coth r_2 \right) (\e^3)^2
\right] .
\end{eqnarray*}
 Plugging these expressions into (\ref{hmm}),
we see that the components of $\schou$  with
respect to the dual frame $\{ \e_j\}$ of $\{ \e^j\}$ satisfy
\begin{eqnarray*}
\schou_{11} &=& (\alpha+1) \sin^2 \varphi + \beta \cos^2 \varphi
 \\&>&\sin^2 \varphi + \beta ,\\
\schou_{22} &=& (\alpha - \beta  )\cos^2 \varphi -\sin^2 \varphi
\\&>&-\sin^2 \varphi , \\
\schou_{33} &=&  (\alpha+1) - (\beta + 1 ) \cos^2 \varphi
 \\&>&(\beta + 1)\sin^2 \varphi ,\\
\schou_{34} &=& \schou_{43} \\&=&
\gamma \sin \varphi \cos \varphi,\\
\schou_{44} &=&  \sin^2 \varphi + \beta \cos^2 \varphi ,\\&&\\
\schou_{jk}&=&0 ~ \mbox{  otherwise,}
\end{eqnarray*}
where
$\alpha : = \coth r_1 + \coth r_2 - 2 $,
$\beta := \frac{\coth^2 r_1 + \coth^2 r_2-2}{\coth r_1 + \coth r_2},$
and $\gamma := \coth r_1 - \coth r_2$ satisfy
 $\alpha  > \beta > |\gamma |$.

Now since
 $\schou_{33}$ and $\schou_{44}$ both exceed $\sin^2\varphi$, and since
\begin{eqnarray*}
 \left|
\begin{array}{cc}
\schou_{33}-\sin^2\varphi&\schou_{34}\\
\schou_{43}&\schou_{44}-\sin^2\varphi\end{array}
\right| &>& (\beta \sin^2 \varphi)
(\beta \cos^2 \varphi)-\gamma^2\sin^2 \varphi
\cos^2 \varphi\\
&=& (\beta^2-\gamma^2) \sin^2 \varphi
\cos^2 \varphi
\\&\geq&0,
\end{eqnarray*}
the eigenvalues $[\schou_{jk}]$ in the
 $\e_3\e_4$-plane  exceed $\sin^2\varphi$. Hence
three of the eigenvalues of $[\schou_{jk}]$ exceed
$\sin^2\varphi$, whereas the remaining eigenvalue $\schou_{22}$
is greater than $-\sin^2\varphi$. The sum of the lowest
two eigenvalues of  $\schou$, calculated with respect to
any metric in the fixed conformal class,
 is therefore
positive   on the  domain of our moving frame. But since this
domain is actually
dense,  it follows that  the Ricci operator is
non-negative on the entirety
 of $\cpt\approx 2{\bf CP}_2$.

Since $\schou_{11} > (\schou_{33}+
\schou_{44})/2 = (\alpha/2) +  \sin^2\varphi$,
the largest two eigenvalues of $[\schou_{jk}]$
are at least $(\alpha/2) +\sin^2\varphi$
on the domain of our frame,
and the sum of the lowest and  third lowest
 eigenvalues of $[\schou_{jk}]$
 therefore exceeds
 $\alpha/2$ on this region.
However, the frame  $\{ \e_j\}$ we have been using is only
{\em conformally} orthonormal with respect to
$g=e^{2f}V(h+V^{-2}\theta^2)$. We now remedy this
by introducing the $g$-orthonormal frame
$\e_j':=e^{-f}V^{-1/2}\e_j$,  with respect to which
 the components of
$\schou$ become
$$\schou_{jk}'=e^{-2f}V^{-1}\schou_{jk}=
\frac{2e^{r_1+r_2}}{\coth r_1 +\coth r_2}\schou_{jk}.$$
If $\mu_1\leq \mu_2\leq \mu_3\leq\mu_4$ are the eigenvalues of
$[\schou_{jk}']$,  we therefore have
\begin{eqnarray*}
\mu_1  + \mu_3 &>& \frac{\alpha e^{r_1+r_2}}{\coth r_1 +\coth r_2}
= e^{r_1+r_2} \frac{e^{2r_1}+e^{2r_2}-2}{e^{2(r_1+r_2)}-1}\\
&\geq& e^{r_1+r_2} \frac{2e^{r_1+r_2}-2}{e^{2(r_1+r_2)}-1}
=\frac{2}{1+e^{-(r_1+r_2)}} > 1.
\end{eqnarray*}
Because  the domain of our frame is dense,
the continuity of the spectrum  therefore
implies that the sum $\mu_1+\mu_3$ of the lowest and
third lowest eigenvalues of $\schou$,
calculated with respect to $g$,  is at least $1$
 on all of $\cpt$.
 The sum $\mu_1+\mu_2$ of the two lowest eigenvalues of $\schou$
can thus vanish only at points at which
 $\schou$ does  not have an eigenvalue of multiplicity 3,
and the Ricci curvature of
$g$ is  therefore strongly positive on all of $\cpt$.

\begin{cor}
Any self-dual metric of positive scalar curvature on ${\bf CP}_2\#
{\bf CP}_2$ is conformal to a metric of strongly positive Ricci
curvature and non-negative Ricci operator.
\end{cor}

\begin{proof} Any   self-dual
conformal class  on ${\bf CP}_2 \# {\bf CP}_2$
 with a representative of
positive scalar curvature is  \cite[p. 251]{leb2}
 of the form $C_{p_1, p_2}$.
\end{proof}

With this success in hand, it seems reasonable, more generally,
to  investigate the
Ricci curvature of  metrics of the form $e^{2f}(Vh+V^{-1}\theta^2)$
on $n{\bf CP}_2$, where $V$ is defined by
(\ref{defV}) and
$$
f = - \frac{r_1 + \cdots + r_n}{n}.
$$
In fact, a rough picture is not difficult to obtain   when the
   points $p_1, \ldots , p_n\in {\cal H}^3$  are extremely close
together. Indeed,
 consider a sequence of configurations of $n$
distinct points in ${\cal H}^3$  which converges to the degenerate
configuration consisting of a
single point $p\in {\cal H}^3$ counted with multiplicity
$n$. On the complement of any  ball about $p$, the
curvature of these metrics will converge uniformally to
that of the orbifold metric corresponding to
$V=1+nG$ and $f=-r$, where   $r$ is the hyperbolic
distance from $p$ and $G=(\coth r - 1)/2$.
  But (\ref{eq12}) predicts that the
  Ricci tensor  of this orbifold limit is
$$\Ric = \zeta \left[ dr^2 + (V^{-1}\theta)^2\right]
+\eta ~ (h-dr^2),$$
where
\begin{eqnarray*}
\zeta &=&
  \frac{\coth r - 1}{2 + n(\coth r -1)} (4 + 3n\coth r - n),
\\
\eta &=&  \frac{\coth r - 1}{2 + n(\coth r -1)} (8 + 3n\coth r - 5n).
\end{eqnarray*}
Observe that
$\eta$ is positive everywhere on ${\cal H}^3$ iff $n \leq 4$, and that
 $\lim_{r\to \infty}\eta/\zeta = 0$ if $n=4$; moreover, we always have
$\zeta \geq \eta$.
Hence the Ricci curvature of this
orbifold limit   is everywhere positive
if and only if $n \leq 3$.
(When $n=4$, it is still non-negative, but fails to be positive
along $\widehat{S}^2$.)
In short, the only encouraging news pertains
to the $n=3$ case, where   the above computation
will help us to prove the following:

\begin{thm}
 If $p_1, p_2,  p_3\in {\cal H}^3$ are
 nearly geodesically collinear and are sufficiently close to each other,
then the conformal class
$C_{p_1, p_2, p_3}$
of self-dual metrics on
$3{\bf CP}_2$
contains a metric with positive Ricci curvature.\label{3cp2}
\end{thm}

To produce self-dual metrics with positive
Ricci curvature on  $3{\bf CP}_2$, we start
 with the above singular model and pull the centers
$p_1, p_2, p_3$ slightly apart, keeping them geodesically collinear.
Outside a neighborhood of $p$, the Ricci curvature remains positive
by our previous computation.
   Theorem \ref{3cp2} is thus implied by   the following:

\begin{Lemma}\label{le4.2}
There exists an $\varepsilon > 0$ such that, for
all collinear configurations $\{ p_1, p_2, p_3\}
\subset {\cal H}^3$, the  Ricci
curvature of $g$ is positive on the inverse image of
$\cup_{j=1}^3 B_{\varepsilon}(p_j)$.
\end{Lemma}

\begin{proof}
Ignoring bounded terms,
$
Ddf \sim -\frac{1}{3} \sum_j \frac{1}{r_j} \left(h - d{r_j}^2\right)
$,
$
\laph f \sim -\frac{2}{3}
\left(\frac{1}{r_1} + \frac{1}{r_2} + \frac{1}{r_3}\right)
$,
$
V \sim \frac{1}{2}
\left(\frac{1}{r_1} + \frac{1}{r_2} + \frac{1}{r_3}\right)
$, and
$
dV \sim -\frac{1}{2}
\left( \frac{dr_1}{{r_1}^2} + \frac{dr_2}{{r_2}^2}
 + \frac{dr_3}{{r_3}^2} \right)
$.
Equation (\ref{eq12}) therefore tells us that
\begin{eqnarray*}
6V \Ric &\sim& \left[
2\left(\frac{1}{r_1} + \frac{1}{r_2} + \frac{1}{r_3}\right)^2
- \left\langle
\frac{dr_1}{{r_1}^2} + \frac{dr_2}{{r_2}^2} + \frac{dr_3}{{r_3}^2},
dr_1 + dr_2 + dr_3 \right\rangle \right] h\\
&& + 2\left(\frac{1}{r_1} + \frac{1}{r_2} + \frac{1}{r_3}\right)
\sum_j \frac{1}{r_j} \left(h - d{r_j}^2\right)
+
2\left( \frac{dr_1}{{r_1}^2} + \frac{dr_2}{{r_2}^2}
 + \frac{dr_3}{{r_3}^2} \right)
\odot (dr_1 + dr_2 + dr_3)\\
&& + \left[
2\left(\frac{1}{r_1} + \frac{1}{r_2} + \frac{1}{r_3}\right)^2
+ \left\langle
\frac{dr_1}{{r_1}^2} + \frac{dr_2}{{r_2}^2} + \frac{dr_3}{{r_3}^2},
dr_1 + dr_2 + dr_3 \right\rangle \right] (V^{-1} \theta)^2\\
&& - 2 \star \left[
\left( \frac{dr_1}{{r_1}^2} + \frac{dr_2}{{r_2}^2}
 + \frac{dr_3}{{r_3}^2} \right)
\wedge (dr_1 + dr_2 + dr_3)\right] \odot V^{-1}\theta,
\end{eqnarray*}
where $\sim$ means that the difference between the
left- and right-hand sides is  of order $g_0=Vh+V^{-1}\theta^2$
 on $\cup_{j=1}^3 B_{\varepsilon}(p_j)$.
Letting $\widehat{R}$
denote the right-hand side of the above expression, it will thus
suffice for us to show that $\widehat{R}$
dominates $Vg_0=V^2h+\theta^2$, since $\Ric$ will then dominate
$\frac{1-CV^{-1}}{6}g_0$ for some constant $C$,
and so will be positive-definite on $\cup_{j=1}^3 B_{\varepsilon}(p_j)$
for  $\varepsilon$  sufficiently small.

Because we are only considering collinear configurations,
$dr_1+dr_2+dr_3 \neq 0$ on ${\cal H}^3\setminus \{ p_1 , p_2, p_3\}$,
and we may let $\e^1$ be the unit covector in this direction.
At any given point, choose $\e^2$ so that the
$dr_j$ are all linear combinations
of $\e^1$ and $\e^2$:
$$dr_j= \cos \varphi_j  ~ \e^1 + \sin \varphi_j ~ \e^2.$$
Extend this to an oriented orthonormal coframe $\{ \e^1, \e^2, \e^3\}$ for
$h$, and set $\e^4=V^{-1}\theta$. Then, letting
$\kappa := \sum_j \cos \varphi_j$,
  the components of $\widehat{R}$
with respect to the  dual frame $\{ \e_j\}$ are
\begin{eqnarray*}\widehat{R}_{11} &=&
\sum_j \frac{1}{{r_j}^2} \left( 2 + \kappa\cos\varphi_j + 2\sin^2\varphi_j
\right)
+ 2\sum_{j<k}\frac{1}{r_jr_k}\left(2+\sin^2\varphi_j + \sin^2\varphi_k\right),
\\
\widehat{R}_{22} &=&
\sum_j \frac{1}{{r_j}^2} \left( 2 - \kappa\cos\varphi_j + 2\cos^2\varphi_j
\right)
+ 2\sum_{j<k}\frac{1}{r_jr_k}\left(2+\cos^2\varphi_j + \cos^2\varphi_k\right),
\\
\widehat{R}_{12} &=& \widehat{R}_{21}\\ &=&
 \sum_j  \frac{1}{{r_j}^2} (\kappa - 2\cos\varphi_j)\sin\varphi_j
- 2\sum_{j<k} \frac{1}{r_jr_k}
(\cos\varphi_j\sin\varphi_j + \cos\varphi_k\sin\varphi_k),
\\
\widehat{R}_{33} &=& \sum_j \frac{1}{{r_j}^2} (4-\kappa\cos\varphi_j)
+ 8 \sum_{j<k}\frac{1}{r_jr_k},
\\
\widehat{R}_{44}&=& \sum_j \frac{1}{{r_j}^2} (2+\kappa\cos\varphi_j)
+ 4 \sum_{j<k}\frac{1}{r_jr_k},
\\
\widehat{R}_{34} &=& \widehat{R}_{43}=
 \sum_j \frac{1}{{r_j}^2} \kappa\sin\varphi_j,
\\&&\\
\widehat{R}_{jk} &=& 0 ~ \mbox{ otherwise.}
\end{eqnarray*}

We now just need to show that the eigenvalues of $[\widehat{R}_{jk}]$ are
all bigger than $V^2$. To do this,   first
notice  that $\sum_k\sin\varphi_k = 0$, and so
\begin{eqnarray*}
\kappa\cos(\varphi_j-2\vartheta) &=&
\left(\sum_k\cos\varphi_k\right) \cos(\varphi_j-2\vartheta) -
\left(\sum_k\sin\varphi_k\right) \sin(\varphi_j-2\vartheta)\\
&=& \sum_k\cos(\varphi_j+\varphi_k-2\vartheta)
\end{eqnarray*}
for any $\vartheta$.
Thus
 $$\cos^2\vartheta\widehat{R}_{11} + 2\cos\vartheta\sin\vartheta
 \widehat{R}_{12}
+\sin^2\vartheta\widehat{R}_{22}=\sum_j\frac{a_j(\vartheta)}{{r_j}^2}
+\sum_{j<k}\frac{a_{jk}(\vartheta)}{r_jr_k}$$
and
$$\cos^2\vartheta\widehat{R}_{33} + 2\cos\vartheta\sin\vartheta
 \widehat{R}_{34}
+\sin^2\vartheta\widehat{R}_{44}=\sum_j\frac{b_j(\vartheta)}{{r_j}^2}
+\sum_{j<k}\frac{b_{jk}(\vartheta)}{r_jr_k},$$
where
\begin{eqnarray*} a_j(\vartheta)
&:=&\cos^2\vartheta \left( 2 + \kappa\cos\varphi_j
+ 2\sin^2\varphi_j\right)
+\sin^2\vartheta \left( 2 - \kappa\cos\varphi_j + 2\cos^2\varphi_j\right)\\
&&
+2\cos\vartheta\sin\vartheta (\kappa-2\cos\varphi_j)\sin\varphi_j\\
&=& 3 + \kappa\cos(\varphi_j-2\vartheta) - \cos(2\varphi_j-2\vartheta)\\
&=& 3 + \sum_{k\neq j} \cos(\varphi_j+\varphi_k-2\vartheta)\\
&\geq& 1,
\\b_{j}(\vartheta)
&:=&\cos^2\vartheta (4-\kappa\cos\varphi_j)
+2\cos\vartheta\sin\vartheta\kappa\sin\varphi_j
+ \sin^2\vartheta(2+\kappa\cos\varphi_j)\\
&=& 3 + \cos 2\vartheta - \kappa\cos(\varphi_j+2\vartheta)\\
&=& 3 - \sum_{k\neq j}\cos(\varphi_j-\varphi_k+2\vartheta)\\
&\geq& 1,
\\a_{jk}(\vartheta)
&:=&2\cos^2\vartheta\left(2+\sin^2\varphi_j + \sin^2\varphi_k\right)
+2\sin^2\vartheta\left(2+\cos^2\varphi_j + \cos^2\varphi_k\right)\\
&&-4\sin\vartheta\cos\vartheta(\cos\varphi_j\sin\varphi_j
 + \cos\varphi_k\sin\varphi_k)\\
&=& 6 - \cos(2\varphi_j-2\vartheta) - \cos(2\varphi_k-2\vartheta)\\
&\geq& 4>2,
\\b_{jk}(\vartheta)
&:=&
8\cos^2\vartheta + 4\sin^2\vartheta = 4 + 4\cos^2\vartheta
\\&\geq& 4>2.
\end{eqnarray*}
Hence every eigenvalue of $[\widehat{R}_{jk}]$ exceeds
$\sum_j\frac{1}{{r_j}^2}+\sum_{j<k}\frac{2}{r_jr_k}=
\left(\sum_j\frac{1}{r_j}\right)^2$, and hence
exceeds $V^2$ on $\cup_j B_{\varepsilon}(p_j)$
for any $\varepsilon < \frac{1}{2}$.  The result follows. \end{proof}

\end{document}